\documentclass[11pt]{article}
\usepackage[a4paper,margin=1in]{geometry}
\usepackage{amsfonts,amssymb,amsmath}
\usepackage{graphicx}
\usepackage{psfrag}
\usepackage{calc}
\usepackage{epic}

\def\R{ {\mathbb{R} }}
\def\M{ {\mathcal{M} }}

\newtheorem{theorem}{Theorem}
\newtheorem{proposition}{Proposition}

\newtheorem{conjecture}{Conjecture}

\def\loynes{\theta^*}

\newlength{\noteWidth}
\long\def\notes#1{\ifinner
             {\tiny #1}
             \else
              \marginpar{\parbox[t]{\noteWidth}{\raggedright\tiny #1}}
               \fi}

\def\notes#1{\typeout{#1 !!!}}  

\begin{document}
\author{
Ken R. Duffy\thanks{Hamilton
Institute, National University of Ireland Maynooth, 
Ireland. 
E-mail: \texttt{ken.duffy@nuim.ie}}
\and 
Sean P. Meyn\thanks{Dept. of Electrical and Computer Engineering,
and the Coordinated Science Laboratory, University of Illinois
at Urbana-Champaign, Urbana, IL, 61801, U.S.A.
E-mail: \texttt{meyn@illinois.edu}}
}
\title{Estimating Loynes' exponent}
\date{$28^{\rm th}$ October 2009}
\maketitle

\begin{abstract}
Loynes' distribution, which characterizes the one dimensional
marginal of the stationary solution to Lindley's recursion, possesses
an ultimately exponential tail for a large class of increment
processes. If one can observe increments but does not know their
probabilistic properties, what are the statistical limits of
estimating the tail exponent of Loynes' distribution? We conjecture
that in broad generality a consistent sequence of non-parametric
estimators can be constructed that satisfies a large deviation
principle.  We present rigorous support for this conjecture under
restrictive assumptions and simulation evidence indicating why we
believe it to be true in greater generality.

\end{abstract}

\section{Introduction}
If $\{X(n)\}$ is a stationary, ergodic process with $E(X(1))<0$,
then Loynes \cite{Loynes62} proved that there is a stationary process
satisfying Lindley's recursion \cite{Lindley52}, 
\begin{align*}
W(n+1) = [W(n)+X(n+1)]^+ \text{ for all } n\in\mathbb{Z},
\end{align*}
and that all other solutions couple to it in almost surely finite
time. A one dimensional marginal of the stationary solution is equal
to Loynes' distribution, the distribution of the random variable
$W = \sup_{n\geq 0} \sum_{i=1}^n X(i)$, where the empty sum
$\sum_{i=1}^0 X(i)$ is defined to be  $0$.

We are interested in estimating the tail behavior of Loynes'
distribution. Consider the partial sums process $\{S(n)/n\}$, where
$S(n):=X(1)+\cdots+X(n)$. Its associated scaled Cumulant Generating
Function (sCGF) $\lambda$ and associated Loynes' exponent $\loynes$
are defined as follows,
\notes{K.D. For Loynes' exponent, if we replace your proposed
$\theta^\circ$ with $\loynes$, we are consistent with Glynn \& Whitt.
Okay?}
\begin{align}
\label{def:scgf}
	\lambda(\theta) := 
	\lim_{n\to\infty}\frac1n \log E\bigl[e^{\theta S(n)}\bigr] \, ,
\qquad 
\loynes:= \sup(\theta:\lambda(\theta)\leq0).
\end{align}
In an advance on the generality of earlier results, Glynn and Whitt
\cite[Theorem 1]{Glynn94} proved the following, which justifies
the terminology Loynes' \textit{exponent}.
\begin{theorem}[\cite{Glynn94}]
\label{thm:glynnwhitt}
Assume that: (i) $\{X(n)\}$ is strictly stationary; (ii) $\lambda(\theta)$
is finite in a neighborhood of $\loynes$, differentiable at $\loynes$
with $\lambda(\loynes)=0$ and $\lambda'(\loynes)>0$; and (iii) that
$E(\exp(\loynes S(n))<\infty$ for all $n\geq1$. Then Loynes' distribution
has an ultimately exponential tail
\begin{align}
\label{eq:tailass}
        \lim_{n\to\infty} \frac 1n \log P(W>n) =-\loynes \, .
\end{align}
\end{theorem}
Note that Lelarge \cite[Proposition 1]{Lelarge08} has recently shown
that condition (ii) can be relaxed. The uniform bound, condition
(iii), ensures that large $W$ is caused by the cumulative behavior
of a collection of increments $\{X(n)\}$ rather than a single
exceedingly large increment.

Theorem \ref{thm:glynnwhitt} is of practical interest as it says
that, for large $n$ and a sizeable collection of increment processes
$\{X(n)\}$, $P(W>n) \sim e^{-n\loynes}$. If $\loynes$ is known,
this provides an estimate of the likelihood of long waiting times
or large queue-lengths, and hence quality-of-service metrics. This
has generated interest in estimating $\loynes$ on-the-fly from
observations to predict queueing behavior in ATM networks. For
example, see early work of Courcoubetis et al. \cite{Courcoubetis95}
and Duffield et al. \cite{Duffield95A}.

\paragraph{Acknowledgment}
{\footnotesize
Financial support from the National Science Foundation (NSF CCF
07-29031), AFOSR (FA9550-09-1-0190) and Science Foundation Ireland
(07/IN.1/1901) is gratefully acknowledged.  Any opinions, findings,
and conclusions or recommendations expressed in this material are
those of the authors and do not necessarily reflect the views of
NSF, AFOSR or SFI.}
\notes{K.D. Added SFI grant}

\section{Estimating Loynes' exponent, a conjecture}

If you can observe a sequence of consecutive increments,
$X(1),\ldots,X(n)$ (or, alternatively, $W(1),\ldots,W(n)$) and wish
to estimate Loynes' exponent, $\loynes$, what are the statistical
properties that your estimator can have?

\begin{conjecture}
\label{con:1}
In broad generality (i.e.\ conditions similar to those in Theorem
\ref{thm:absurd}) without knowing anything further about the process
$\{X(n)\}$, one can build a sequence of estimates $\{\loynes(n)\}$
that satisfy a Large Deviation Principle (LDP) \cite{Dembo98}. That
is, for all Borel sets
\begin{align}
\label{eq:deltan}
-\inf_{x\in B^\circ} J(x) \le
\liminf_{n\to\infty}\frac1n\log P(\loynes(n)\in B) \le
\limsup_{n\to\infty}\frac1n\log P(\loynes(n)\in B) \le
-\inf_{x\in \bar{B}} J(x),
\end{align}
where $B^\circ$ denotes the interior of $B$ and $\bar{B}$ denotes
the closure of $B$, and $J:[0,\infty]\mapsto[0,\infty]$ is a good
rate function (lower semi-continuous and with compact level sets).
Moreover, one can construct consistent estimates ($J(x)=0$
if and only if $x=\loynes$).
\end{conjecture}

Equation \eqref{eq:deltan} can be compared with \cite[Proposition~11.3.4]{CTCN}
  or
\cite[Theorem 3]{Duffy09}.  The following is a corollary of the latter.  
\notes{S.M. 1. We did need assumptions, and I felt I couldn't package
this up without 'proposition'.  2. Is this right? I just attempted
to derive the sqrt bound without success!\\
K.D. You're right on both points. Regarding the latter, if we
restrict the $w$ range then the result is as stated in your book 
and as deduced by your recent scaling arguments.}
\begin{proposition}[\cite{Duffy09}]
\label{t:DM09}
Suppose that the sequence $\{X(n)\}$ is i.i.d. with $P(X(1)>0)>0$
and assume that $\lambda(\theta) =  \log E\bigl[e^{\theta X(1)}]$
is finite in a neighborhood of the origin. Then with
\[
S_w(n):=W(1)+\cdots+W(n),
\]   
the process $\{S_w(n)/n^2\}$ satisfies the LDP with a good rate function
$J$.  That is, 
\[
P\Bigl\{\sum_{i=1}^n W(i) > n^2x\Bigr\}\sim \exp(-nJ(x))
\]
with $J(0)=0$, 
$J(x)>0$ for all $x>0$, 
and $J(x)<\infty$ for some
$x>0$.
\end{proposition}

Consequently,  we are conjecturing that estimating the tail exponent
of Loynes' distribution is easier than estimating its mean.

{\bf Estimation schemes and rigorous evidence in support of the conjecture.}
A range of approaches to estimating $\loynes$ can be deduced from
the literature. In increasing order of directness, examples include
the following. (I) For the GI/G/1 queue, techniques have been
developed for estimating the whole of Loynes' distribution
\cite{Pitts94} from which $\loynes$ can be deduced. (II) Estimators
of the sCGF have been proposed that use either a frequentist approach
if $\{X(n)\}$ satisfies a mixing condition \cite{Duffield95A}, or
a Bayesian approach if $\{X(n)\}$ is i.i.d. \cite{Ganesh98}\cite{Ganesh00}
or a Markov chain \cite{Paschalidis01}. From these, estimates of
$\loynes$ can be obtained (c.f. \cite{Hall92}). (III) Direct extremal
estimators based on observations of the $\{W(n)\}$, such as
$\log(n)/\max(1,W(1),\ldots,W(n))$, have been studied
\cite{Berger95}\cite{Zeevi04}.

Our limited rigorous evidence for the validity of Conjecture
\ref{con:1} is based on two frequentist sCGF estimation schemes.
We show that the conjecture holds if the non-overlapping partial
sums of the increments $\{X(n)\}$ are i.i.d. and bounded for a fixed
block size. This result is deduced from \cite{Duffy05}.  Moving
away from independence, we also present a new result: if $\{X(n)\}$
forms a finite state irreducible Markov chain, then the conjecture
holds. For the conjecture to be established in generality, new ideas
are needed to extend to unbounded and non-independent increments
that take more than finite values, but - arguably - the unboundedness
is more technically challenging. We hope this will become apparent
in the exposition that follows.

The first result is based on the estimation scheme proposed in
\cite{Duffield95A}.  For an integer $B<\infty$, construct the blocked
process
\begin{equation}
Y(i)= \sum_{j=(i-1)B+1}^{iB} X(j). 
\label{e:Duffield95A}
\end{equation}
Select $B$ sufficiently large that you believe
the blocked process $\{Y(n)\}$ is close to being i.i.d. If the
process $\{Y(n)\}$ was i.i.d., then $\lambda(\theta)$ in equation
\eqref{def:scgf} reduces to $\lambda(\theta) = B^{-1}\log E(\exp(\theta
Y(1)))$.  Given observations $X(1),\ldots,X(n)$, this suggests using
the MLE for $\lambda(\theta)$\footnote[2]{This is known to be a biased
estimator \cite{Ganesh96}.} and $\loynes$: 
\begin{align*}
	\hat\lambda(n,\theta) := \frac 1B \log
		\left(\frac{1}{\lfloor n/B\rfloor} \sum_{i=1}^{\lfloor
		n/B\rfloor} e^{\theta Y(i)}\right)
\text{ and finally }
	\loynes(n) := \sup\{\theta: \hat\lambda(n,\theta)\leq0\}.
\end{align*} 
A central limit theorem for $\{\loynes(n)\}$ is proved in
\cite{Duffield95A}. As a corollary to a result regarding a related
estimation problem, \cite[Theorem 2]{Duffy05} proves that the
sequence of estimates $\{\loynes(n)\}$ satisfy the LDP under 
less restrictive conditions than those of the following theorem,
but does not establish its consistency.

\begin{theorem}[\cite{Duffy05}]
\label{thm:absurd}
If, for some $B$, $\{Y(i)\}$ is i.i.d. and $Y(i)$ takes values in a
closed, bounded subset $\Sigma$ of the real line that does not
include an open ball around the origin, then Conjecture \ref{con:1}
holds.
\end{theorem}
A sketch of the proof of Theorem \ref{thm:absurd} is as follows.
By Sanov's Theorem, the empirical laws $\{L(n)\}$, defined by
$L(n):=n^{-1}\sum_{i=1}^n 1_{Y(i)}$ for $n\geq1$, satisfy the LDP
in $\M_1(\Sigma)$, the space of probability measures on $\Sigma$,
equipped with the topology of weak convergence. If $Y(1)$ has measure
$\mu$, then the good rate function for the LDP is the relative
entropy $H(\nu|\mu)$.
As $\Sigma$ is bounded, for each $\theta\in\R$ the function 
$x\mapsto\exp(\theta x)$ is continuous and bounded. Thus if
$\mu_n\implies \mu$ in $\M_1(\Sigma)$, then 
$\log\mu_n(\exp(\theta x))\to\log\mu(\exp(\theta x))$ in $\R$.
As point-wise convergence of convex functions implies uniform
convergence on bounded subsets, by the contraction principle
$\{\hat\lambda(n,\cdot)\}$ satisfies the LDP in the space of convex
functions equipped with the topology of uniform convergence on
bounded subsets. To obtain the LDP for the tail exponent estimates
$\{\loynes(n)\}$, one considers the continuity of
$\hat\lambda(n,\cdot)\mapsto\sup(\theta:\hat\lambda(n,\theta)\leq0)$ and
applies Puhalskii's extension of the contraction principle
\cite[Theorem 2.1]{Puhalskii95}. As $H(\nu|\mu)=0$ if and only if
$\nu=\mu$, consistency follows from the variational form of the
rate function for $\{\loynes(n)\}$ that is given by the contraction
principle.

The need to exclude $Y(n)$ taking values in an open ball around the
origin is due to an artifact of the estimation scheme. If
$Y(1)=\cdots=Y(n)=0$, then $\hat\lambda(n,\theta)=0$ for all $\theta$
and $\loynes(n)=\infty$. However, in the topologically nearby situation
where $Y(1)=Y(2)=\cdots=Y(n)=\epsilon>0$, $\loynes(n)=0$. That is,
the estimation scheme possesses a discontinuity and the assumption
is imposed to avoid it.

Our second result is based on a frequentist version of the estimator
used in \cite{Paschalidis01}. Consider a finite state Markov
chain $\{X(n)\}$ with an irreducible transition matrix
$\Pi=(\pi_{i,j})\in[0,1]^{M\times M}$ and taking values
$\{f(1),\ldots,f(M)\}$.
For each $\theta\in\R$, define the matrix $\Pi_\theta= \Pi D_\theta$,
where $D_\theta$ is the matrix with diagonal entries 
$\exp(\theta f(1)), \exp(\theta f(2)),\ldots, \exp(\theta f(M))$
and all off-diagonal entries equal to zero. The sCGF of the partial
sums process $\{S(n)/n\}$ can be identified as
$\lambda(\theta)=\log\rho(\Pi_\theta)$, where $\rho$ is the spectral
radius \cite[Theorem 3.1.2]{Dembo98}. With $0/0:=0$, this suggests
that one constructs the MLE for $\Pi$, $\hat\Pi(n)$, defined by
\begin{align*}
\hat\pi(n)_{i,j} :=  
	\left(\sum_{k=1}^n 1_{\{(X(k-1),X(k))=(i,j)\}} \right) /
		\left(\sum_{k=1}^n 1_{\{X(k-1)=i\}}\right),
\end{align*}
and then estimates $\lambda(\theta)$ and $\loynes$ by
$\hat\lambda(n,\theta) = \log\rho\left(\hat\Pi(n)_\theta\right)$ 
and
$\loynes(n) = \sup(\theta:\hat\lambda(n,\theta)\leq0)$.
\begin{theorem}
\label{thm:absurd1}
If $\{X(n)\}$ is a finite state Markov chain with an irreducible
transition matrix $\Pi$ and $f(i)\neq0$ for all $i\in\{1,\ldots,M\}$,
then Conjecture \ref{con:1} holds true.
\end{theorem}
Theorem \ref{thm:absurd1} follows the arguments of Theorem
\ref{thm:absurd} once we establish that $\{\hat\lambda(n,\cdot)\}$
satisfies the LDP. By \cite[Theorem 3.1.13]{Dembo98}, the empirical
laws of the transitions $\{L_2(n)\}$,
$L_2(n):=n^{-1}\sum_{i=1}^n1_{(X(n-1),X(n))}$, satisfy the LDP in
$\{(i,j):\pi_{i,j}>0\}$. With $\phi=(\phi_1,\ldots,\phi_M)$ being
the stationary distribution of $\Pi$, the rate function defined by
\cite[equation (3.1.14)]{Dembo98} is zero only at the values
$\phi_i\pi_{i,j}$. If $\pi_{i,j}=0$, then $\hat\pi(n)_{i,j}=0$.
For all $(i,j)$ such that $\pi_{i,j}>0$, we have that 
$\hat\pi(n)_{(i,j)}$ can be expressed as a ratio of integrals
against $L_2(n)$:
\begin{align*}
\hat\pi(n)_{(i,j)} = 
	\left(L_2(n)(1_{(x,y)=(i,j)})\right)/
 	\left(L_2(n)(1_{(x,y)=(i,\cdot)})\right).
\end{align*}
Thus the estimate $\hat\Pi(n)$ is a continuous construction from
$L_2(n)$ so that contraction principle can be applied and the
estimates $\{\hat\Pi(n)\}$ satisfy the LDP with a rate function
$H:[0,1]^{M\times M}\mapsto[0,\infty]$ that satisfies $H(A)=0$ if
and only if $A=\Pi$.  If the sequence of matrices $A_n$ converge
entry-wise to $A$, then by the continuous dependence of eigenvalues
on entries \cite{Horn91},
$\log\rho(A_nD_\theta)$ converges to
$\log\rho(AD_\theta)$ point-wise for all $\theta$ and therefore
uniformly on compact subsets of $\theta$. Thus, from the contraction
principle, $\{\hat\lambda(n,\cdot)\}$ satisfies the LDP and the
rest of the proof follows as for Theorem \ref{thm:absurd}.

Following the logic in \cite[Theorem 1]{Duffy05}, as an aside we
note that based on either of these two estimation schemes one can
construct estimates of the rate function for the partial
sums process $\{S(n)/n\}$. Defining 
$\hat I(n,x):=\sup_{\theta}(\theta x-\hat\lambda(n, \theta))$, 
if $\{\hat\lambda(n,\cdot)\}$ satisfies the LDP, then it
can be shown that $\{\hat I(n,\cdot)\}$ satisfies the LDP in the
space of $\R\cup\{\infty\}$ valued convex functions equipped with
the Attouch-Wets topology \cite{Attouch83}\cite{Attouch86}. That
is, there is a LDP for estimating large deviation rate functions
and, under the conditions of Theorems \ref{thm:absurd} and
\ref{thm:absurd1}, the estimates are consistent.

\begin{figure}
\begin{center}
\includegraphics[width=\textwidth*\real{0.45}]{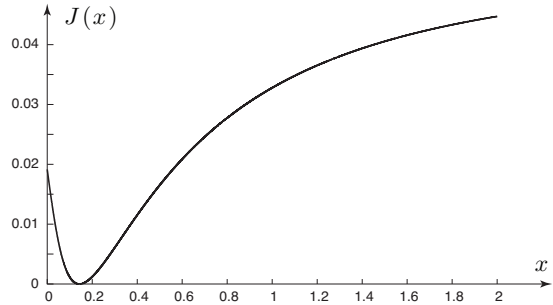}
\caption{Rate function for $\{\loynes(n)\}$ based on frequentist
Markov chain estimator. Markov increments $\{X(n)\}$ on $\{-1,+1\}$
with $\alpha=1/16$, $\beta=3/16$ and $\loynes=\log(15/13)$}
\label{fig:2state}
\end{center}
\end{figure}

{\bf Example.} Assume that the increments process $\{X(n)\}$ 
forms a two-state Markov chain on the state space $\{-1,+1\}$
with transition matrix
\begin{align*}
\Pi = \left( 
	\begin{array}{cc}
	1-\alpha & \alpha\\
	\beta & 1-\beta 
	\end{array}
	\right)\, .
	\qquad
\text{Then,  }
	\phi = \left( \frac{\beta}{\alpha+\beta}, 
		\frac{\alpha}{\alpha+\beta}
		\right)
\text{ and }
	\loynes = \log\left(\frac{1-\alpha}{1-\beta}\right),
\end{align*} 
where $0<\alpha < \beta<1$. 
\notes{S.M. "contraction from" might be confusing\\ 
K.D Agreed. Is this better?}
The sequence $\{\hat\Pi(n)\}$ satisfies the LDP with a good rate
function $H$. It can be deduced from \cite[equation (3.1.14)]{Dembo98}
and the variational expression given by the contraction principle
that $H$ is finite only at matrices of the form
\begin{align*}
A = \left( 
	\begin{array}{cc}
	1-a & a \\
	b & 1-b 
	\end{array}
	\right),
\text{ where } a,b\in(0,1)
\end{align*}
in which case 
\begin{align*}
H(A)
	 &=
	\frac{b}{a+b} \left(
		(1-a)\log\left(\frac{1-a}{1-\alpha}\right)
		+a\log\left(\frac{a}{\alpha}\right)
		\right) 
	 +\frac{a}{a+b} \left(
		b\log\left(\frac{b}{\beta}\right)
		+(1-b)\log\left(\frac{1-b}{1-\beta}\right)
		\right).
\end{align*}
Consequently, the rate function for $\{\loynes(n)\}$ is given by
the one dimensional optimization
\begin{align*}
J(x) = \inf\left(H(A(a,b):\log\left(\frac{1-a}{1-b}\right)=x\right) 
	= \inf_{a\in(0,1)}H(A(a,1-(1-a)e^{-x})).
\end{align*}
While $J(x)$ cannot be determined in closed form, it can be readily
calculated numerically. Figure~\ref{fig:2state} provides an example
for given parameters; its non-convex nature for large $x$ is
apparent.
\notes{K.D. Typo corrected: first $\delta$ is $\beta =P(X(n+1)=-1|X(n)=+1)$
and second $\delta$ is Loynes' exponent $\loynes$.}  

{\bf Simulation evidence to support the conjecture.}
For the conjecture to be substantiated, the conditions under which
Theorems \ref{thm:absurd} and \ref{thm:absurd1} hold need to be
significantly extended to cope with more general dependence structure
of the increments and, perhaps more challengingly, to remove the
boundedness assumptions.  A prototypical example of the later is
where Lindley's recursion describes the waiting times at the $D/M/1$
queue. That is, the increments $\{X(n)\}$ are i.i.d. with
$P(X(1)>x+1/\beta) = \exp(-\alpha x)$ for $x\geq-\beta^{-1}$ and
$\alpha>\beta$.  This example satisfies the conditions of Theorem
\ref{thm:glynnwhitt} with $\loynes(<\alpha)$ being the positive
solution of the transcendental equation
\begin{equation}
\log\left(\alpha/(\alpha-\loynes)\right)-\loynes/\beta=0,
\label{e:DM1loynes}
\end{equation}
which can
be readily solved numerically. 

\begin{figure}
\begin{center}
\includegraphics[width=\textwidth*\real{.95}]{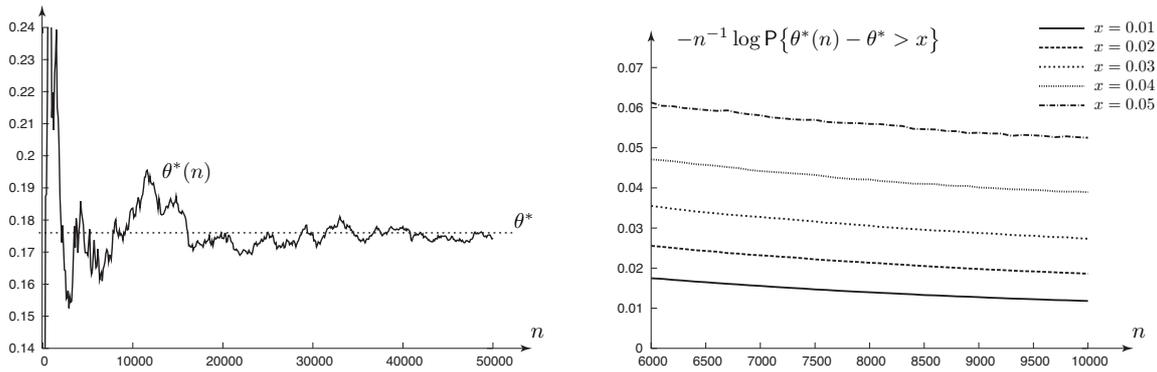}
\caption{D/M/1 queue with $\alpha=1$ and $\beta=10/11$, giving
$\loynes\approx0.176$.  Shown on the left hand side are estimates of $\loynes$ 
based on equation \eqref{e:Duffield95A} with $B=1$.
Plots shown on the right hand side are estimates of the rate function  for 
this estimator evaluated at $\loynes +x$ for several values of $x$.}
\label{fig:dm1}
\end{center}
\end{figure}

This example is delicate because the tail of increments decay
exponentially, albeit with a larger exponent than $\loynes$.  Even
though this example breaks the hypotheses of Theorem \ref{thm:absurd},
we can implement the estimation scheme in \cite{Duffield95A} and
see how it performs.  Fixing $B=1$, for a single realization
$X(1),\ldots,X(50,000)$ the left hand side of Figure~\ref{fig:dm1}
plots $\loynes(n)$ as a function of $n$ as well as the actual Loynes' exponent 
$\loynes$ solving
\eqref{e:DM1loynes}.
This plot indicates that the estimates are converging to the correct
value. The right-hand plot is an attempt to consider the existence
of an LDP. For $2\times 10^5$ independent simulations and a range
of values of $x$, it plots $n^{-1}$ times the logarithm of the
number of samples such that $\hat\loynes(n)-\loynes>x$, which we
expect to converge to the rate function for the tail exponent
estimates.  This figure is suggestive of the existence of an LDP
despite the departure from the boundedness conditions in Theorems
\ref{thm:absurd} and \ref{thm:absurd1}.

\notes{S.M. Is my caption ok?  
\\
When we get the final figure I will change the y-axis to read
$-\frac{1}{n} P\bigl\{ \hat\theta^\circ(n)-\theta^\circ>x \bigr\}$\\
K.D. Great! New figure attached, but sims still running!}

{\bf Further questions.} Conjecture \ref{con:1} is challenging, but
Theorem \ref{thm:glynnwhitt} has been extended considerably so that
more difficult questions can be asked. In a development of the sCGF
approach in \cite{Glynn94}, \cite{Duffield95} considers the case
where the partial sums $\{S(n)/n\}$ satisfy the LDP at a non-linear
speed.  By reconsidering the general scaling problem in terms of
rate functions rather than sCGFs \cite[Theorem 2.2]{Duffy03} extends
the results further while also correcting a lacuna (the omission
of an assumption in the vein of Theorem \ref{thm:glynnwhitt} (iii)).

\begin{theorem}[\cite{Duffy03}]
Assume that $\{S(n)/n^A\}$ satisfies the LDP at speed $n^V$ with
rate function $I$ (so that, roughly speaking, 
$P(S(n)>x n^A)\sim \exp(-I(x) n^V)$) and define
$\loynes := \inf_{x>0} x^V I\left(1/x^A\right)$.
If, in addition, $n^{-V} \log P(S(n)>xn^A) \leq -\loynes$ for all
$n$ and all $x$ sufficiently large, then
\begin{align*}
        \lim_{n\to\infty} \frac {1}{n^{V/A}} \log P(W>n) =-\loynes.
\end{align*}
\end{theorem}
For example, this Theorem holds if $\{X(n)\}$ is i.i.d. with a
Weibull distribution \cite{Nagaev79} or if $\{S(n)\}$ corresponds
to sampled fractional Brownian motion or a sampled
two state process with Weibull sojourn times \cite{Duffy08}.

Among our questions are: can one simultaneously estimate $V/A$ while
estimating $\loynes$?  What impact does the non-linear scaling have
on the properties of estimators? 

\bibliography{muse}
\bibliographystyle{plain}

\end{document}